\documentclass[12pt]{article}
\usepackage{color}
\usepackage{amsmath,amsfonts,amsthm}
\usepackage{algorithmic}
\newcommand {\eq} [1] {\begin{equation}\label{#1}}
\newcommand {\en} {\end{equation}}
\newcommand {\cA}       {{\cal A}}
\newcommand {\cB}       {{\cal B}}
\newcommand {\cC}       {{\cal C}}
\newcommand {\cD}       {{\cal D}}
\newcommand {\cE}       {{\cal E}}

\newcommand {\cK}       {{\cal K}}
\newcommand {\cL}       {{\cal L}}

\newcommand {\cP}       {{\cal P}}
\newcommand {\cQ}       {{\cal Q}}
\newcommand {\cR}       {{\cal R}}
\newcommand {\cS}       {{\cal S}}
\newcommand {\cT}       {{\cal T}}
\newcommand {\cU}       {{\cal U}}
\newcommand {\cV}       {{\cal V}}
\newcommand {\cW}       {{\cal W}}

\newcommand {\cZ}       {{\cal Z}}
\newcommand {\eproof}
      {\space
        {\ \vbox{\hrule\hbox{\vrule height1.3ex\hskip0.8ex\vrule}\hrule}}
        \par}
\newcommand {\C}        {{\mathbb C}}
\newcommand {\R}        {{\mathbb R}}

\newcommand {\Rn}       {\R^n}
\newcommand {\Rm}       {\R^m}

\newcommand {\Rnn}      {\R^{n,n}}
\newcommand {\Rnm}      {\R^{n,m}}
\newcommand {\Rmm}      {\R^{m,m}}
\newcommand {\Rmn}      {\R^{m,n}}

\newcommand {\mat}      [1] {\left[\begin{array}{#1}}
\newcommand {\rix}          {\end{array}\right]}
\newtheorem{theorem}           {Theorem}
\newtheorem{lemma}    [theorem]{Lemma}

\newtheorem{remark}            {Remark}

\newtheorem{algorithm}         {Algorithm}
\newcommand {\rank}     {\mathop{\rm rank}\nolimits}

 \catcode`@=11
 \font\tenex=cmex10 % math extension
 \newdimen\p@renwd
 \setbox0=\hbox{\tenex B} \p@renwd=\wd0 % width of the big left [
 \def\bmat#1{\begingroup \m@th
   \setbox\z@\vbox{\def\cr{\crcr\noalign{\kern2\p@\global\let\cr\endline}}%
     \ialign{$##$\hfil\kern2\p@\kern\p@renwd&\thinspace\hfil$##$\hfil
       &&\quad\hfil$##$\hfil\crcr
       \omit\strut\hfil\crcr\noalign{\kern-\baselineskip}%
       #1\crcr\omit\strut\cr}}%
   \setbox\tw@\vbox{\unvcopy\z@\global\setbox\@ne\lastbox}%
   \setbox\tw@\hbox{\unhbox\@ne\unskip\global\setbox\@ne\lastbox}%
   \setbox\tw@\hbox{$\kern\wd\@ne\kern-\p@renwd\left[\kern-\wd\@ne
     \global\setbox\@ne\vbox{\box\@ne\kern2\p@}%
     \vcenter{\kern-\ht\@ne\unvbox\z@\kern-\baselineskip}\,\right]$}%
   \null\;\vbox{\kern\ht\@ne\box\tw@}\endgroup}
 \catcode`@=12
\begin{document}
\title{{Regularization of Port-Hamiltonian Descriptor Systems}}
    \author{ Chu Delin \\*
                 Department of Mathematics\\
                National University of Singapore \\*
                 Singapore 119076\\*
                Email:  matchudl@nus.edu.sg.
    \and
Volker Mehrmann \\*
      Institut f\"{u}r Mathematik\\*
      MA 4-5, TU Berlin \\*
      Str. des 17. Juni 136, D-10623 Berlin, FRG \\*
      Email: mehrmann@math.tu-berlin.de.}

\maketitle
\begin{abstract}
We study the regularization problem for  port-Hamiltonian descriptor systems %$E\dot x =A x+B u$, $y=Cx $
by proportional and/or derivative output feedback. Necessary and sufficient conditions are given, which guarantee that there exist  output feedbacks such that the  closed-loop system is regular, has index at most one, and is still port-Hamiltonian with desired rank properties.
{All results are derived based on condensed forms which may be computed in a numerically reliable way using orthogonal transformations.}
\end{abstract}



\section{Introduction}\label{intro}
Port-Hamiltonian systems, see \cite{RasCSS20,Sch09,SchJ14}, provide an ideal  model class for the automated modelling, control and simulation of real world physical systems.  The physical properties are encoded in the structure of the mathematical model, see \cite{RasCSS20,SchJ14}. These physical properties of the systems
include balance and conservation laws, the shaping of energy-storage and energy-dissipation, as well as the interpretation of controller systems as virtual system components. When algebraic constraints are present in the physical system, then it is more appropriate to use the class of port-Hamiltonian des\-criptor systems, see  \cite{MehS23} for the most general definition in the linear time invariant case, \cite{BeaMXZ18} for linear time-invariant systems and \cite{MehM19,Sch13} for nonlinear systems, and a survey in \cite{MehU23}. Due to automated modeling procedures that include redundant equations, the resulting
port-Hamiltonian descriptor systems may not be regular, i.e. for a given input the solutions to initial value problems  may not exist or are not unique.  Furthermore the system may be of index higher than one, i.e. the solution depends on derivatives of the input function. Then impulses may arise in the response of the system if the control
is not sufficiently smooth.

How to deal with non-regular or high index problems is well understood for general descriptor systems. If this is possible, then one   uses feedback to make the system regular and  of index at most one.  This topic has been studied extensively in many different versions. For general linear time-invariant systems the existence of regularizing state or output feedbacks is discussed in \cite{OzcL90,ShaZ87}. The construction of these feedbacks with orthogonal transformations is discussed in \cite{BunMN92,BunMN94,ByeGM97,ChuCH98,ChuH99,NicC15} and a survey is given in \cite{BunBMN99}.  For linear time-varying descriptor systems a summary is presented in \cite{KunM24} and for general nonlinear descriptor systems in
\cite{CamKM12,KunM24}.
Various applications, see e.g. \cite{AltMU21,BanS21,Sch09}, strongly motivate us to study the regularization problem for port-Hamiltonian descriptor systems, since existing regularization procedures in general do not preserve the port-Hamiltonian structure.
Hence, it is important from a  theoretical and practical point of view to develop theory and numerical reliable methods for the regularization and index reduction preserving the port-Hamiltonian structure. For a special class of linear time-invariant port-Hamiltonian descriptor systems and using only proportional feedback, this research topic has recently been studied  in \cite{ChuM25}. In this paper we discuss this topic for general linear time-invariant port-Hamiltonian descriptor systems and we also include derivative output feedback.

The paper is organized as follows. The regularization problem via output feedback for
 port-Hamiltonian descriptor systems is introduced in Section~\ref{S2} as well as some preliminaries. The main results are presented in Section~\ref{S3} and concluding remarks are given in Section~\ref{S4}.

\section{Preliminaries and Problem Statement}\label{S2}
We  study port-Hamiltonian descriptor systems of the form
{
\begin{eqnarray}
E {\dot x} &=& Ax+Bu,  \label{1.1}  \\
         y &=&  Cx,    \label{1.1o}
\end{eqnarray}
}
in a time interval $\mathbb I\subset \mathbb R$,
where $E, A\in \Rnn$, $B\in \Rnm$ and $C\in \Rmn$,   $x: \mathbb I\to \Rn$ is the state, $y: \mathbb I \to \Rm$ is the output,  $u:\mathbb I \to \Rm$ is the input or control of the system. {
For a system of the form \eqref{1.1}-\eqref{1.1o} to be port-Hamiltonian, the coefficient matrices have to satisfy,
\begin{equation}\label{phrep}
\mat{cc} A & B \\ C & 0 \rix =\mat{cc} (J-R)Q & G-P \\ (G+P)^TQ & 0 \rix,
\end{equation}
with $E, J, R, Q \in \Rnn, G, P\in \Rnm$,  $J=-J^T$, %\Gamma:=\mat{cc} J & G \\ -G^T & 0 \rix=-\Gamma^T$,
and the further relations
\begin{equation}
\label{positive} Q^TE=E^TQ\geq 0, \quad Q^TRQ=Q^TR^TQ\geq 0, \ \ Q^TP=0,
\end{equation}
see \cite{MehU23}.

The quadratic function $ \mathcal H(x)={1\over 2} x^T E^TQx$ is called the \emph{Hamiltonian} which can  be interpreted as the energy stored in system. It satisfies the \emph{power balance equation}, see \cite{MehM19,MehU23},
  \begin{equation}\label{eq:powerBalanceEq}
    \frac{d}{dt}\mathcal H(x) = -\mat{c} x\\ u \rix^T\mat{cc} Q^T R Q & Q^TP\\
P^T Q &0 \rix\mat{c} x\\ u\rix + y^Tu
  \end{equation}
along any solution $x$ and for any input $u$. Here the first term on the right hand side is the dissipated energy and the second term is the supplied energy.

\begin{remark} \label{rem:qtp=0}{\rm
It is well-known that  port-Hamiltonian descriptor systems satisfy the dissipation inequality
\begin{equation}\label{eq:dissinq}
    \frac{d}{dt}\mathcal H(x) =  y^Tu,
  \end{equation}
which is referred to as passivity in control. For this to hold in the case of systems without feedthrough term, it is necessary that
\begin{equation}\label{dissinq}
\mat{cc} Q^T R Q & Q^TP\\
P^T Q &0 \rix\geq 0
\end{equation}
which implies that $Q^TP=0$ and $Q^TRQ\geq 0$. This implies in particular, that the uncontrolled system (with $u=0$) is semidissipative.
}
\end{remark}

Note that $E$ and $Q$ are allowed to be singular. Note that the present methods are also directly applicable for problems  with feedthrough term, which can always be formulated as \eqref{1.1}, see \cite{MehU23}. However, we consider only problems without feedthrough term.

}
The response of the descriptor system (\ref{1.1}) can be described in terms of the eigenstructure of the \emph{matrix pencil}
$\,\alpha E - \beta A\,,$ which we abbreviate by $(E,A)$. The system is called \emph{regular}
if the pencil $\,(E, A)\,$ is regular, i.e.,
\begin{equation} \label{1.3}
    {\rm det}(\alpha E-\beta A)\not= 0 {\ \ \mbox {for some}}  \ \ (\alpha, \beta)\in \C^2.
\end{equation}
The \emph{generalized eigenvalues }of a regular pencil $(E, A)$ are defined by the pairs $\,(\alpha_j\,,\beta_j)\,\in\C^2\backslash\{0,0\}\,$ such that
\begin{equation} \label{1.3b}
    {\rm det}(\alpha_j E-\beta_j A) = 0\,, \ \ \ \ j = 1, 2,\ldots,n\,.
\end{equation}
If $\,\beta_j \not= 0\,,$ the eigenvalue pairs are said to be \emph{finite} with value given by $\,\lambda_j = \alpha_j / \beta_j\,$ and otherwise, if $\,\beta_j = 0\,,$ then the pair is said to be an \emph{infinite} eigenvalue.
The maximum number of finite eigenvalues that a pencil $(E, A)$ can have is less than or equal to the rank of $\,E\,.$

In this paper, we denote a full column rank matrix with its columns spanning the right  nullspace of the matrix $M$ by $\cS_\infty(M)$ and with
its columns spanning the left nullspace of $M$ by  $\cT_\infty(M)$, respectively.

If the system (\ref{1.1}) is regular, then for given consistent initial values, {i.e., initial values that satisfy the algebraic equations present in the system,} the existence and uniqueness of classical solutions to the dynamical equations is guaranteed, see \cite{KunM24}.
In the regular case, the solutions can be characterized in terms of the Kronecker Canonical Form (KCF) which states that there exist nonsingular matrices
$X$ and $Y$ (representing right and left generalized eigenvectors and generalized eigenvectors of the system pencil, respectively)  such that
\begin{equation}\label{KCF} X(sE-A)Y=\mat{cc} sI-J & 0 \\ 0 & sN-I \rix,
\end{equation}
where the eigenvalues of  $J$ coincide with the finite eigenvalues of the pencil and $N$ is a nilpotent matrix
such that for $i>0$ it holds that $N^i=0$, $N^{i-1}\not=0$. the matrix $N$ is corresponding to the infinite eigenvalues. The \emph{index} of a descriptor system,
denoted by ${\rm \mbox{\rm ind }}(E, A)$, is defined to be the degree $i$ of nilpotency of the matrix $N$, i.e., the index of the system is the dimension of the largest block associated with an infinite eigenvalue in the KCF (\ref{KCF}).
Note that the system  (\ref{1.1}) is regular and has index at most one if and only if it has exactly $ {\rm  rank}(E)$ finite eigenvalues.

Due to the port-Hamiltonian structure, the index of a  regular
%completely controllable and observable
port-Ha\-mil\-to\-ni\-an descriptor system of the form (\ref{1.1}) is at most $2$  \cite{MehMW18,MehS23,MehU23},  but its index may indeed be $2$.

% Therefore, the main purpose of this paper is to  establish the complete theory and to develop numerically reliable methods for  the regularization of port-Hamiltonian descriptor systems via derivative and proportional output feedback.

As a special case  for port-Hamiltonian descriptor systems with $Q=I$, the identity,
%of the form
%
%\begin{eqnarray*}
%   E\dot x &=& (J-R)x+Bu, \\
%   y &=& B^Hx,
%\end{eqnarray*}
%
the regularization problem by proportional output feedback has been studied recently in \cite{ChuM25}.   In this paper, we extend these results to the more general
form of port-Hamiltonian descriptor systems as in (\ref{1.1}).

We consider proportional and derivative output feedbacks of the form
\[
u=Fy-K\dot y+v,
\]
where $F, K$ are feedback matrices so that the resulting system
has the form
\begin{eqnarray} \label{f2}
   (E+BKC)\dot x &=& (A+BFC)x+Bv, \nonumber\\
              y &=& Cx,
\end{eqnarray}
with desired properties.  Here, proportional feedback control is achieved with $K=0$ and
derivative feedback control corresponds to the case $F=0$.
%Derivative and proportional state feedback control corresponds to the
%case
%\begin{equation} \label{state}  u(t)=Fx(t)-K\dot x(t)+v(t).
%\end{equation}
Proportional feedback {changes} the system matrix $A$, while derivative feedback alters the  matrix $E$. Therefore, different properties of the system can be achieved using different
feedback combinations.

We study the following problems:

a)  \emph{Regularization problem for descriptor system (\ref{1.1})   by {proportional} output feedback:}
Determine a matrix $F$  such that
 the  pencil $(E, A+BFC)$ is regular, $\mbox{\rm ind }(E, A+BFC)\leq 1$,  and the resulting system
\begin{eqnarray}
   E \dot x(t) &=& (A+BFC)x(t)+Bv(t), \label{sys1}\\
              y(t) &=& Cx(t),  \nonumber
\end{eqnarray}
is  still  port-Hamiltonian, {i.e., satisfies \eqref{phrep} with the properties \eqref{positive}.}

b)  \emph{Regularization problem for descriptor system (\ref{1.1}) by derivative output feedback:}
Determine a matrix  $K$  such that
 the  pencil $(E+BKC, A)$ is regular, $\mbox{\rm ind }(E+BKC, A)\leq 1$, and the resulting system
\begin{eqnarray}
   (E+BKC) \dot x(t) &=& Ax(t)+Bv(t), \label{sys2}\\
              y(t) &=& Cx(t),  \nonumber
\end{eqnarray}
is  still  port-Hamiltonian, {i.e., satisfies \eqref{phrep} with the properties \eqref{positive}.}

c)  \emph{Regularization problem of descriptor system (\ref{1.1})  by derivative and proportional output feedback:}
Determine matrices $F$ and $K$ such that   the  pencil $(E+BKC, A+BFC)$ is regular, $\mbox{\rm ind }(E+BKC, A+BFC)\leq 1$, and the resulting system
\begin{eqnarray}
 (E+BKC) \dot x(t) &=& (A+BFC)x(t)+Bv(t), \label{sys3}\\
              y(t) &=& Cx(t),    \nonumber
\end{eqnarray}
is  still  port-Hamiltonian, {i.e., satisfies \eqref{phrep} with the properties \eqref{positive}.}

%These problems above will be studied in the next section.

\section{Main Results}\label{S3}

One of the main difficulities for solving the regularization problem of port-Hamiltonian descriptor systems is to preserve the port-Hamiltonian structure  by feedbacks.

The following results present necessary and sufficient solvability conditiond for the regularization of the port-Hamiltonian descriptor system (\ref{1.1}) by proportional output feedback and derivative output feedback,
respectively.

{The following result presents the answer to Problem a).}

\begin{theorem}\label{theorem-1}
Let system (\ref{1.1}) be port-Hamiltonian. Then there exists a matrix $F$ such that the pencil $(E, A+BFC)$ is regular, $\mbox{\rm ind }(E, A+BFC)\leq 1$, and the closed-loop system (\ref{sys1}) is port-Hamiltionan
 if and only if
\begin{equation}\label{con-1}
\rank \mat{cc} E & A\cS_\infty(E) \\ 0 & C\cS_\infty(E) \rix=n. % \rank \mat{ccc} E & A \cS_\infty(E) & B \rix=n
\end{equation}
\end{theorem}

\proof Condition (\ref{con-1}) follows from the necessary conditions for the  regularization problem of general descriptor system (\ref{1.1})   by output feedback without the port-Hamiltonian requirement, see {Theorem 3.3 in \cite{ChuCH98}}. Hence, the necessity follows.

To show the sufficiency, we give a constructive proof, that can be implemented as numerical method, to show
that under the condition (\ref{con-1}) there exists a matrix $F\in \Rmm$ such that $(E, A+BFC)$ is regular and of index at most one and the closed-loop system (\ref{sys1}) is port-Hamiltonian.

We first compute orthogonal matrices $U$ and $V$  such that
\begin{equation}\label{EABC-0}  U EV=\bmat{& n_1 & n_2 & n_3 \cr
         n_1    & E_{11} & 0 & 0 \cr
         n_2    & 0 & 0 & 0          \cr
        n_3    & 0 & 0 & 0 \cr},  \
    U AV =\bmat{& n_1 & n_2 & n_3 \cr
                  n_1 &  A_{11} & A_{12} & A_{13} \cr
                  n_2 &  A_{21} & A_{22} & 0 \cr
                  n_3  &  A_{31} & 0 & 0 \cr},
\end{equation}
where $E_{11}$ and $A_{22}$ are nonsingular.  It is well-known, see e.g. \cite{ChuM25}, that for a  port-Hamiltonian descriptor system of the form \eqref{1.1} the matrix $Q$ in (\ref{positive}) satisfies
\begin{equation}\label{inequality}  E^TQ=Q^TE\geq 0, \quad  -A^TQ-Q^TA\geq 0, \quad C-B^TQ=0.
\end{equation}
{

If $Q$ is partitioned analogously, then the first inequality in  \eqref{inequality} implies that  $Q_{12}=0$
and $Q_{13}=0$, since $E_{11}$ is nonsingular. Then from the second inequality it follows that $A_{22} Q_{23}=0$ and hence $Q_{23}=0$, since $A_{22}$ is nonsingular so that we have}
\[
UQV=\bmat{& n_1 & n_2 & n_3 \cr
               n_1 & Q_{11} & 0 & 0 \cr
              n_2  & Q_{21} & Q_{22} & 0 \cr
              n_3 & Q_{31} & Q_{32} & Q_{33} \cr}.
\]
Partition
\[
UB=\bmat{ &       \cr
              n_1 &  B_1 \cr
              n_2 & B_2 \cr
              n_3 & B_3 \cr},
              \quad CV=\bmat{& n_1 & n_2 & n_3 \cr
            & C_1 & C_2 & C_3\cr}
\]
accordingly.
Using the condition (\ref{con-1}) it follows that $C_3$ is of full column rank, so we can compute an orthogonal matrix $W$ such that
\begin{equation} W^T \mat{ccc} C_1 & C_2 & C_3 \rix
                   =\bmat{& n_1 & n_2 & n_3 \cr
                   m-n_3   & C_{11} & C_{12} & 0 \cr
                       n_3   &  C_{21} & C_{22} & C_{23} \cr}, \label{C1}
\end{equation}
where $C_{23}$ is nonsingular because $\rank(C_{23})=\rank(C_3)$. {This follows from condition \eqref{con-1}, since multiplication with $\mathcal S_\infty(E)$ projects onto the last two block columns. }

Partition
\[
\mat{c} B_1 \\ B_2 \\ B_3 \rix W=\bmat{& m-n_3 & n_3 \cr
         n_1 & B_{11} & B_{12} \cr
        n_2 & B_{21} & B_{22} \cr
        n_3 & B_{31} & B_{32}  \cr}.
\]
Then the condition $C=B^TQ$  implies that
\[
B_{32}^TQ_{33}=C_{23}, \ B_{31}^TQ_{33}=0,
\]
which, together with the nonsingularity of $C_{23}$, gives that $B_{32}$ and $Q_{33}$ are nonsingular and $B_{31}=0$.

Then for
\[
F=-W\mat{cc} 0 & 0  \\ 0  & F_{22} \rix W^T
\]
with
\[
F_{22}\in \R^{n_3\times n_3}, \quad  F_{22}=F_{22}^T>0,
\]
a direct calculation yields that $(E, A+BFC)$ is regular and of index at most one. Furthermore,
\[
A+BFC=(J-R)Q+BFB^TQ=[J-(R-BFB^T)]Q,
\]
and
\[
Q^T(R-BFB^T)Q=Q^T(R-BFB^T)^TQ= Q^TRQ-Q^TBFB^TQ\geq 0.
\]
Hence, the closed-loop system (\ref{sys1}) is port-Hamiltonian.
\eproof

For derivative feedback we have an analogous result {which presents the answer to Problem b).}
\begin{theorem}\label{theorem-2}
Let system (\ref{1.1}) be port-Hamiltonian. Then there exists a matrix $K$ such that the pencil $(E+BKC, A)$ is regular, $\mbox{\rm ind }(E+BKC, A)\leq 1$,
$\rank(E+BKC)=\max_{\hat K\in\Rmm} \rank(E+B\hat KC)$, and the closed-loop system (\ref{sys2}) is port-Hamiltonian
 if and only if
\begin{equation}\label{con-2}
\rank \mat{cc} E & A\cS \left (\mat{c} E \\ C \rix \right ) \\ C & 0 \rix=\rank \mat{ccc} E & A\cS_\infty\left (\mat{c} E \\ C \rix \right ) & B \rix=n.
\end{equation}
\end{theorem}
\proof
Condition (\ref{con-2}) follows from the necessary conditions for the  regularization problem of general descriptor systems (\ref{1.1})   by derivative output feedback without port-Hamiltonian requirement,
{see Theorem 3.1 in \cite{ChuCH98}}. Hence, the necessity follows.

For the port-Hamiltonian descriptor system  (\ref{1.1}),  it holds that
\[
\rank \mat{c} E \\ C \rix\leq \rank \mat{cc} E & B \rix.
\]
{This follows since $C=B^TQ$ and $E^TQ=Q^TE\geq 0$.}
Thus,
\[
\max_{\hat K\in\Rmm} \rank(E+B\hat KC)=\rank \mat{c} E \\ C \rix.
\]
To show the sufficiency, we {present a construction of such a regularizing feedback that can be implemented as numerical method. For this, we
first  separate  the part of the system that cannot be influenced by derivative output feedback, and that is already of index at most one, } by determining orthogonal matrices $U$ and $V$,
such that
\[
UEV=\bmat{ & n_1 & n_2 \cr
                \hat n_1 & E_{11} & 0 \cr
                \hat n_2 & E_{21} & E_{22} \cr}, \quad UAV=\bmat{& n_1 & n_2 \cr
                \hat n_1 & A_{11} & 0 \cr
                \hat n_2 & A_{21} & A_{22} \cr}, \quad
      CV=\bmat{ & n_1 & n_2 \cr
                       & C_1 & 0 \cr},
\]
where
\[
\rank \mat{c} E_{11} \\ C_1 \rix=n_1,
\]
and the index at most one system
\[
\rank(sE_{22}-A_{22})=\hat n_2, \ \mbox{for all}\  s\in \C.
\]
{This can be achieved by first forming a full column rank decomposition of $C$ and then  to use the staircase algorithm in the kernel to separate the infinite Jordan blocks of size one, see \cite{DemK93}.}
Partition
\[
UB=\bmat{&   \cr
               \hat n_1 &  B_1 \cr
             \hat n_2  &  B_2 \cr}
\]
accordingly. Then condition (\ref{con-2}) is equivalent to
\[
n_2=\hat n_2, \quad  E_{22}=0, \quad \rank(A_{22})=n_2, \ \rank \mat{cc} E_{11} & B_1 \rix=n_1.
\]
Thus, condition (\ref{con-2}) implies that
\[
\rank \mat{c} E \\ C \rix=n_1.
\]
Because $\rank\mat{c} E_{11} \\ C_1 \rix=\rank \mat{cc} E_{11} & B_1 \rix=n_1$,  it is easy to determine a matrix $K$ { (for example by taking a scaled identity of size $n_1$)} such that
\[
K=K^T\geq 0, \quad \rank(E_{11}+B_1KC_1)=n_1.
\]
Obviously, $(E+BKC, A)$ is regular,  $\mbox{\rm ind }(E+BKC, A)\leq 1$, and $\rank(E+BKC)=\rank \mat{c} E \\ C \rix$. Moreover, with $Q$ as in (\ref{inequality}), then
\begin{eqnarray*}
(E+BKC)^TQ&=&E^TQ+C^TK^TB^TQ=E^TQ+Q^TBKB^TQ\\
   &=&Q^T(E+BKC)\geq 0.
\end{eqnarray*}
Hence,  the closed-loop system (\ref{sys2}) is port-Hamiltionan.
\eproof

The necessary and sufficient conditions in Theorems \ref{theorem-1} and \ref{theorem-2} are easily verified and the desired feedback matrices $F$ and $K$ can be computed easily by applying orthogonal transformations which can be implemented in a numerically stable.

Note that the necessary and sufficient conditions for the  regularization of descriptor system (\ref{1.1})   by proportional output feedback,  presented in \cite{ChuCH98},
without the port-Hamiltonian requirement, are the condition (\ref{con-1}) and
\begin{equation}\label{con-1-1}
\rank \mat{ccc} E & A \cS_\infty(E) & B \rix=n.
\end{equation}
Because of the port-Hamiltonian structure, the condition (\ref{con-1}) implies the condition (\ref{con-1-1}), so the condition (\ref{con-1-1}) is not needed for the
port-Hamiltonian descriptor system (\ref{1.1}).

In the following construction,  we assume without loss of generality that $C$ in (\ref{1.1}) is of full row rank,  i.e., $\rank(C)=m$.

In order to derive solvability conditions for the  regularization of port-Hamiltonian descriptor system (\ref{1.1}) by derivative and proportional output feedback, we need the following two lemmas.

\begin{lemma}\label{lemma-3}
 Let the system (\ref{1.1}) be port-Hamiltonian and $\rank\mat{c} E \\ C\rix=n$.  Then there exist
orthogonal matrices $U$, $V$ and $W$ such that
\begin{eqnarray}
\nonumber
U E V &=&\bmat{  & n-r_b & r_e+r_b-n & n-r_e \cr
                     n-r_b   &  E_{11} & E_{12} & 0  \cr
                     r_e+r_b-n   & E_{21} & E_{22} & 0  \cr
                            n-r_e   &  0 & 0 & 0   \cr},  \\
U B W &=&\bmat{ & r_e+r_b-n & n-r_e \cr
                     n-r_b   &   0 & B_{12}  \cr
                     r_e+r_b-n   &  B_{21} & B_{22} \cr
                     n-r_e   &  0 & B_{32}  \cr}, \label{condensedform1}\\
W^T C V&=&\bmat{ & n-r_b & r_e+r_b-n & n-r_e \cr
                     r_e+r_b-n   &  0 & C_{12} & 0 \cr
                     n-r_e   &  C_{21} & C_{22} & C_{23} \cr},
                     \nonumber
\end{eqnarray}
where $B_{21}$, $C_{12}$, $B_{32}$, $C_{23}$, and $E_{11}$ are nonsingular.
\end{lemma}

\proof The proof is given in the appendix.\eproof

Now consider  orthogonal matrices $\cU$ and $\cV$ with partitioning
\[
\cU=\bmat{ & n-r_b & r_e+r_b-n \cr
          n-r_b    & \cU_{11} & \cU_{12} \cr
     r_e+r_b-n   & \cU_{21} & \cU_{22} \cr}, \quad
\cV=\bmat{ & n-r_b & r_e+r_b-n \cr
               n-r_b  & \cV_{11} & \cV_{12} \cr
       r_e+r_b-n   & \cV_{21} & \cV_{22} \cr},
\]
that satisfy
\[
U \mat{c} E_{11} \\ E_{21} \rix=\mat{c} \cE_{11} \\ 0 \rix, \quad \rank(\cE_{11})=n-r_b,
\]
and
\[
\mat{cc} E_{11} & E_{12} \rix \cV=\mat{cc} \tilde \cE_{11} & 0 \rix, \quad \rank(\tilde \cE_{11})=n-r_b.
\]
%$\cU_{22}$ and $\cV_{22}$ are nonsingular \cite{ChuLM03}.
Set
\begin{eqnarray*}
\hat E_{22}&=&(\cU_{21}E_{12}+\cU_{22}E_{22})\cV_{22}, \ \hat B_{21}=\cU_{22}B_{21}, \ \hat B_{22}=\cU_{21}B_{12}+\cU_{22}B_{22},\\
  \hat C_{12}&=&C_{12}\cV_{22}, \quad \hat C_{22}=C_{21}\cV_{12}+C_{22}\cV_{22},
\end{eqnarray*}
and
\[
\mat{ccc} I &     &       \\ \cU_{21} & \cU_{22} & \\ &  & I \rix U A V \mat{ccc} I & \cV_{12} & \\                                                & \cV_{22} &  \\   &  & I \rix
=\bmat{  & n-r_b & r_e+r_b-n & n-r_e \cr
                     n-r_b   &  A_{11} & A_{12} & A_{13}  \cr
                     r_e+r_b-n   & A_{21} & A_{22} & A_{23}  \cr
                     n-r_e   &  A_{31} & A_{32} & A_{33}   \cr}.
\]
Then determine orthogonal matrices $Z$ and $\cZ$ such that
\[
Z \mat{cc} A_{22} & A_{23} \\ A_{32} & A_{33} \rix \cZ=\bmat{ & \mu & r_b-\mu \cr                                                                    \mu           & \cA_{22} & 0 \cr                                           r_b-\mu           & 0            & 0 \cr}, \]
where  $\cA_{22}$ is nonsingular.
\begin{lemma}\label{lemma-4}  Consider a port-Hamiltonian system of the form (\ref{1.1}) and suppose that $\rank \mat{c} E \\ C \rix=n$. Let  $U$, $V$ and $W$ be orthogonal matrices that transform the system to the condensed form \eqref{condensedform1} and let
 $\cW$ be orthogonal  such that
\[
Z \mat{cc} \hat B_{21} & \hat B_{22} \\ 0 & B_{32} \rix \cW=
\bmat{&    \mu      &  r_b-\mu    \cr
\mu     &  \cB_{21} & \cB_{22} \cr
r_b-\mu     &   0           & \cB_{32} \cr},
\]
where  $\cB_{21}$, $\cB_{32}$ are nonsingular.
Then
\[
\cW^T \mat{cc} \hat C_{12} & 0 \\ \hat C_{22} & C_{23} \rix \cZ =
\bmat{ & \mu & r_b-\mu \cr                                                  \mu           & \cC_{12}  & 0 \cr                                          r_b-\mu          &  \cC_{22} & \cC_{23} \cr},
\]
where $\cC_{12}$ and $\cC_{23}$ are nonsingular.
\end{lemma}
\proof The proof is given in the appendix. \eproof

We are now ready to present necessary and sufficient solvability conditions for the regularization  of port-Hamiltonian descriptor systems
by derivative output feedback and derivative
and proportional output feedback (with a given rank of $E+BKC$). {The following result presents the answer to Problem c).}
\begin{theorem}\label{theorem-3} Suppose that the port-Hamiltonian system (\ref{1.1}) is \emph{completely observable},  i.e., $\rank \mat{c} \alpha E-\beta A \\ C \rix=n$
for any $(\alpha,\beta) \in\C^2\backslash \{0,0\}$. {Then following assertions hold.}

(i) There exists a matrix $K$ such that  $(E+BKC, A)$ is  regular,
\[ \mbox{\rm ind }(E+BKC, A)\leq 1,   \ \rank(E+BKC)=r, \]
and the closed-loop system (\ref{sys2}) is port-Hamiltionan if and only if
\begin{equation}\label{R1} n-\rank[ \cT_\infty^T(E\cS_\infty(C))A \cS_\infty (\cT_\infty^T(B)E)] \leq r \leq n,
\end{equation}
and $n-r$ is even when $\rank[ \cT_\infty^T(E\cS_\infty(C))A \cS_\infty (\cT_\infty^T(B)E)]>0$ is even and
\[
(\cT_\infty^T(E\cS_\infty(C))B)^{-1} [\cT_\infty^T(E\cS_\infty(C))A \cS_\infty (\cT_\infty^T(B)E)] (C \cS_\infty \c(T_\infty^T(B)E)^{-1}
\]
is skew-symmetric.

(ii)  There exist matrices $F$ and $K$ such that  $(E+BKC, A+BFC)$ is regular,
\[
\mbox{\rm ind }(E+BKC, A+BFC)\leq 1, \ \rank(E+BKC)=r,
\]
 and the system (\ref{sys3}) is port-Hamiltonian if and only if
\[
n-r_b \leq r \leq n.
\]
\end{theorem}

\proof  Since  (\ref{1.1}) is completely observable, we have  $\rank \mat{c} E \\ C \rix=n$.
Using (\ref{1.7}) we have that
\[
\rank[ \cT_\infty^T(E\cS_\infty(C))A \cS_\infty (\cT_\infty^T(B)E)]=\rank\mat{cc} A_{22} & A_{23} \\                       A_{32} & A_{33} \rix=\mu,
\]
and
\begin{eqnarray*}
&& (\cT_\infty^T(E\cS_\infty(C))B)^{-T} [\cT_\infty^T(E\cS_\infty(C))A \cS_\infty (\cT_\infty^T(B)E)] (C \cS_\infty (\cT_\infty^T(B)E)^{-1} \\
&& =\mat{cc} \hat B_{21} & \hat B_{22} \\ 0 & B_{32} \rix^{-1}\mat{cc} A_{22} & A_{23} \\ A_{32} & A_{33} \rix \mat{cc} \hat C_{12} & 0 \\ \hat C_{22} & C_{23} \rix^{-1}\\
 &&  =\cW \mat{cc} \cB_{21}^{-1} \cA_{22}\cC_{12}^{-1} & 0 \\ 0 & 0 \rix \cW^T.
\end{eqnarray*}
Since $-A^TQ-Q^TA\geq 0$, using (\ref{1.10}), we have  that
\[
\cB_{12}^{-1}\cA_{22}\cC_{12}^{-1}+ (\cB_{12}^{-1}\cA_{22}\cC_{12}^{-1})^T \leq 0.
\]
Hence, there exists an orthogonal matrix $\hat \cP$ such that
\begin{equation}\label{Schur-form}
\hat \cP ( \cB_{12}^{-1}\cA_{22}\cC_{12}^{-1})\hat   \cP^T=\hat \cA_{22},
\end{equation}
is in real Schur form, see \cite{GolV96},  and $\hat \cA_{22}+\hat \cA_{22}^T\leq 0$. Set
\begin{eqnarray*}
\hat X&=&\mat{ccc} I & & \\ & \hat \cP & \\ & & I \rix \mat{ccc} I &  & \\ & \cB_{21} & \cB_{22} \\ & 0 & \cB_{32} \rix^{-1} \mat{cc} I & \\ & Z \rix X, \\
 \hat Y&=&Y \mat{cc} I & \\ & \cZ \rix \mat{ccc} I & & \\ & \cC_{12} & 0 \\  & \cC_{22} & \cC_{23} \rix^{-1} \mat{ccc} I  &  & \\ & \hat \cP^T & \\ & & I \rix.
 \end{eqnarray*}
Then
\[
\hat X Q \hat Y =\mat{ccc} \cQ_{11} &  & \\ & I & \\ & & I \rix,
\]
and for  $K, F\in \R^{m\times m}$ form
\begin{eqnarray*}
W^TKW&=&\mat{cc} \cK_{11} & 0 \\ 0 & 0 \rix+ (\cW \hat \cP) \mat{cc}  K_{11} & K_{12} \\                                                         K_{21} & K_{22} \rix (\cW \hat \cP)^T,  \\ W^TFW=(\cW \hat \cP) \mat{cc}  F_{11} & F_{12} \\                                                         F_{21} & F_{22} \rix (\cW \hat \cP)^T,
\end{eqnarray*}
where $\cK_{11}, K_{11}, F_{11}  \in \R^{(n-r_b)\times (n-r_b)}$, $B_{21}\cK_{11}C_{12}=E_{22}$.
Then we have
\begin{eqnarray*} \hat X  (E+BKC) \hat Y  &=&\mat{ccc} E_{11} & 0 & 0 \\
                                     0 & K_{11} & K_{12} \\ 0 & K_{21} & K_{22} \rix, \ \hat X A\hat Y
= \mat{ccc} A_{11} & \hat \cA_{12} & \hat \cA_{13} \\ \hat \cA_{21} & \hat \cA_{22} & 0 \\ \hat \cA_{31} & 0 & 0 \rix, \\\
 \hat X  (A+BFC) \hat Y  &=&\mat{ccc} A_{11} & \hat \cA_{12}  & \hat \cA_{13} \\
\hat \cA_{21} & \hat \cA_{22}+F_{11}  & F_{12} \\ \hat \cA_{31} & F_{21} & F_{22} \rix.
\end{eqnarray*}

(i) Then $(E+BKC, A)$ is regular, $\mbox{\rm ind }(E+BKC, A)\leq 1$ and $\rank(E+BKC)=r$ if and only if
\[
(\mat{cc} K_{11} & K_{12} \\ K_{21} & K_{22} \rix,
            \mat{cc} \hat \cA_{22} & 0 \\ 0 & 0\rix
\]
is regular, of index at most one and
\[
\rank \mat{cc} K_{11} & K_{12} \\ K_{21} & K_{22} \rix=r-(n-r_b)=r+r_b-n. \]
To show the other parts, observe first that if $(E+BKC, A)$ is regular, $\mbox{\rm ind }(E+BKC, A)\leq 1$ and the closed-loop system (\ref{sys2}) is  port-Hamiltonian, then it follows that
%
%\begin{equation}\label{X}
\[
\mat{cc} K_{11} & K_{12} \\ K_{21} & K_{22} \rix \geq 0,
\]
%\end{equation}
hence the following assertions hold
\begin{eqnarray*} \rank(K_{22}) &=&\rank \mat{cc} K_{21} & K_{22} \rix=r_b-\mu, \\
r&=&n-r_b+\rank \mat{cc} K_{11} & K_{12} \\ K_{21} & K_{22} \rix \geq n-r_b+\rank (K_{22})\\
&=&n-r_b+r_b-\mu
      =n-\mu, \\
K_{21}&=&K_{12}^T, \ K_{11}-K_{12}K_{22}^{-1}K_{12}^T \geq 0, \end{eqnarray*}
and, furthermore, $\cS_\infty(K_{11}-K_{12}K_{22}^{-1}K_{12}^T)^T   \hat  \cA_{22} \cS_\infty(K_{11}-K_{12}K_{22}^{-1}K_{12}^T)$ is nonsingular.

Since any real skew-symmetric matrix with odd dimension is singular, it follows that when $\mu>0$ and $\cA_{22}+\cA_{22}^T=0$, then we have that
$\mu$ is even and $\rank(\cS_\infty(K_{11}-K_{12}K_{22}^{-1}K_{12}^T)^T  \hat \cA_{22} \cS_\infty(K_{11}-K_{12}K_{22}^{-1}K_{12}^T)$
is also even, i.e.,  $n-r$ is even.

Next, let $r$ be any integer satisfying (\ref{R1}), i.e.,  $n-\mu \leq r \leq n$.  We have that $\hat \cA_{22}$ is in real Schur form
\[
\hat \cA_{22}=\mat{cccc} \cT_{1} & *       & * & * \\
                                             & \ddots & * &   *    \\
                                             &           & \cT_k  & *    \\
                                                          &         & & \cD \rix, \]
where $\cD$ is in  real Schur form, $\cD+\cD^T\not=0$, and
\[
\cT_i=\mat{cc} 0 & t_i \\ -t_i & 0 \rix, \quad t_i\not=0, \quad i=1, \cdots, k.
\]
Choose $K_{12}=0$,  $K_{21}=0$, $K_{22}=I$, and let $K_{11}$ be constructed as follows:

 If $2k=\mu> 0$,   then $\cA_{22}+\cA_{22}^T=0$, $r+\mu-n$ is even and if $r+\mu-n=2s$, then set
\[
K_{11}=\mat{cc} I_{2s} &  \\ & 0 \rix.
\]
Otherwise, if $2k<\mu$, then we choose $K_{11}$ as in the following cases:
\begin{eqnarray*}
K_{11}&=&\mat{cc} I_{2s} &    \\ & 0 \rix, \ {\rm when \ } r+\mu-n=2s\leq 2k; \\
K_{11}&=&\mat{ccc} I_{2s} &      &    \\
            & 0    &             \\
            &       & 1      \rix, \ \ {\rm when \ } r+\mu-n=2s+1, \ s\leq k;       \\
K_{11}&=&\mat{ccc} I_{2k} &    &       \\
        & I_{r+\mu-n-2k} &     \\
        &  & 0 \rix, \ {\rm when \ } r+\mu>2k,
\end{eqnarray*}
The the pair $(K_{11}, \hat \cA_{22})$ is regular, $\mbox{\rm ind }(K_{11}, \hat \cA_{22})\leq 1$ and $\rank(K_{11})=r+\mu-n$,
and thus, $(E+BKC, A)$ is regular, $\mbox{\rm ind}(E+BKC, A)\leq 1$ and $\rank(E+BKC)=r$. Moreover,
\begin{equation}\label{P}  \mat{ccc} E_{11} &     &        \\
             & K_{11} &       \\
            &            & I \rix \mat{ccc} \cQ_{11} & &       \\
        &   I   &      \\
        &       & I    \rix
= \mat{ccc} \cQ_{11} &      &       \\
 &   I   &      \\
        &       & I    \rix^T
     \mat{ccc} E_{11} &     &        \\
                & K_{11} &       \\
            &            & I \rix  \geq 0,
\end{equation}
or  equivalently,
\[
(E+BKC)^TQ=Q^T(E+BKC)\geq 0
\]
and, thus, the closed-loop system (\ref{sys2}) is port-Hamiltonian.

(ii) Obviously,  $\rank(E+BGC)\geq r_{eb}-r_b=n-r_b$ for any $K\in \R^{m\times m}$. Let $r$ be an integer with $n-r_b\leq r\leq n$ and choose
$K_{12}=0$, $K_{21}=0$, $K_{22}=I$, $F_{12}=0$, $F_{21}=0$, $F_{11}=-I$, $F_{22}=-I$,
and
\[
F_{11} =\mat{cc} I_{r+\mu-n} & \\ & 0 \rix.
\]
We then have that $(K_{11}, \hat \cA_{22}+F_{11})$ is regular, $\mbox{\rm ind }(K_{11}, \hat \cA_{22}+F_{11})\leq 1$ and $\rank(K_{11})=r+\mu-n$. Hence,
$\rank(E+BKC)=r$, $(E+BKC, A+BFC)$ is regular, and $\mbox{\rm ind }(E+BKC, A+BFC)\leq 1$. In addition, (\ref{P}) holds, and thus
$(E+BKC)^TQ\geq 0$. Moreover,
\[
A+BFC=(J-R)Q+BFB^TQ=[J-(R-BFB^T)]Q,
\]
and
\[
Q^T(R-BFB^T) Q=Q^TRQ-Q^TBFB^TQ \geq 0,
\]
and, therefore, the closed-loop system (\ref{sys3}) is  port-Hamiltonian.
\eproof

The condensed form (\ref{condensedform1}) is  computed using  only orthogonal transformations  and so it can be implemented as a numerically reliable algorithm. The Schur form (\ref{Schur-form}) can be computed as follows:

\begin{algorithm}
i) Compute the QR factorization
\[
\mat{c} \cC_{12} \\ \cA_{22} \rix=\mat{cc} L_{11} & L_{12} \\ L_{21} & L_{22} \rix \mat{c} R \\ 0 \rix,
\]
where $\mat{cc} L_{11} & L_{12} \\ L_{21} & L_{22} \rix$ is orthogonal and  $R$ is nonsingular. Then $L_{11}$ and $L_{22}$ are nonsingular, see e.g. \cite{ChuLM03}, and
\[
\cA_{22}\cC_{12}^{-1}=L_{21}L_{11}^{-1}=-L_{22}^{-T}L_{12}^T, \quad
   \cB_{21}^{-1}\cA_{22}\cC_{12}^{-1}=-\cB_{21}^{-1}(L_{22}^T)^{-1}L_{12}^T.
\]

ii) Compute the QR factorization
\[ \mat{cc} -L_{22}^T\cB_{21} & L_{12}^T\rix =\mat{cc} \cR & 0 \rix \mat{cc} \cL_{11} & \cL_{12} \\ \cL_{21} & \cL_{22} \rix, \]
where $\mat{cc} \cL_{11} & \cL_{12} \\ \cL_{21} & \cL_{22} \rix$ is orthogonal and $\cR$ is nonsingular. Then again $\cL_{11}$ is nonsingular,
see \cite{ChuLM03}, and
\[
\cB_{21}^{-1}\cA_{22}\cC_{12}^{-1}=\cL_{11}^{-1}\cL_{12}.
\]
Compute $\cL_{11}^{-1}\cL_{12}$ by using the Cosine-Sine Decomposition, see \cite{GolV96}, of $\cL_{11}$ and $\cL_{12}$ and then compute the real Schur form (\ref{Schur-form}).
\end{algorithm}
All the computations, and hence, the numerical construction of matrices $K$ and $F$ in Theorem \ref{theorem-2}, can be implemented in a numerically reliable way.

\begin{remark}{\rm
An important feature of Theorems \ref{theorem-1}, \ref{theorem-2} and \ref{theorem-3} and the construction of the feedback matrices
 is that the computation of the matrix $Q$ in (\ref{positive}) is not directly involved, because this would be computationally difficult.}
\end{remark}

\section{Concluding Remarks}\label{S4}
We have presented necessary and sufficient  conditions for the regularization of  port-Hamiltonian descriptor systems by derivative and/or proportional output feedback. All results are derived based on condensed forms that can be implemented using only orthogonal transformations and hence as numerically reliable algorithms.
%Future work will include the study of the stabilization problem of port-Hamiltonian descriptor system (\ref{1.1}). Because the port-Hamiltonian property must be preserved for the closed-loop system, the associated stabilization problem is very difficult and  is still a challenging task.

\section{Appendix}

\subsection {Proof of Lemma~\ref{lemma-3}}

\proof  We again give a constructive proof that {can be directly implemented as numerical method.}
First observe that
\[
m=r_c=\rank(C)=\rank(B^TQ)\leq \rank(B)=r_b\leq m,
\]
so that
\begin{equation}\label{1.5}
r_c=\rank(C)= \rank(B)=r_b=m.
\end{equation}

The construction of  $U$, $V$  and $W$  proceeds as follows:
Determine orthogonal matrices $U_1$ and $V_1$ such that
\[
E^{(1)}=U_1EV_1=\bmat{ & r_e & n-r_e \cr
                          r_e      & E_{11}^{(1)} & 0 \cr
                      n-r_e       & 0                   & 0 \cr}, \quad \rank(E_{11}^{(1)})=r_e.
\]
Partition {analogously}
\[
B^{(1)}=U_1B=\bmat{ &         \cr
                           r_e & B_1^{(1)} \cr
                       n-r_e  & B_2^{(1)} \cr}, \quad C^{(1)}=CV_1=\bmat{ & r_e & n-r_e \cr
            & C_1^{(1)} & C_2^{(1)} \cr}.
\]
Note that by $\rank \mat{c} E \\ C \rix=n$, {$C=B^TQ$} and $E^TQ=Q^TE\geq 0$,  it follows that
\[
\rank(C_2^{(1)})=n-r_e,
\]
and
\[
Q^{(1)}=U_1QV_1=\bmat{ & r_e & n-r_e \cr
                  r_e & Q_{11}^{(1)} & 0 \cr
             n-r_e & Q_{21}^{(1)} & Q_{22}^{(1)} \cr}.
\]
Compute an orthogonal matrix $W$ such that
\[
B_2^{(1)} W=\bmat{ & r_b+r_e-n & n-r_e \cr
                          & 0 &     B_{22}^{(2)} \cr}
\]
and partition
 \[
 B^{(2)}=\mat{c} B_1^{(1)} \\ B_2^{(1)} \rix W=
 \bmat{ & r_b+r_e-n & n-r_e \cr
 r_e              & B_{11}^{(2)} & B_{12}^{(2)} \cr
n-r_e              & 0                   & B_{22}^{(2)} \cr}
\]
accordingly. Then by $C=B^TQ$ it follows
\begin{eqnarray}\label{1.6a}  C^{(2)}&=&W^T \mat{cc} C_1^{(1)} & C_2^{(1)} \rix
 = \bmat{ & r_e & n-r_e \cr
r_b+r_e-n              & C_{11}^{(2)} & 0 \cr
 n-r_e                     & C_{21}^{(2)} & C_{22}^{(2)} \cr}, \\
 C_{22}^{(2)}&=&(B_{22}^{(2)})^TQ_{22}^{(1)}.
\end{eqnarray}
Since $\rank(C_2^{(1)})=n-r_e$, so, $\rank(C_{22}^{(2)})=\rank(C_{2}^{(1)})=n-r_e$ and  consequently  $\rank(B_{22}^{(2)})=n-r_e$ and $Q_{22}^{(1)}$ is nonsingular.
Note that $\rank(C)=\rank(B)=m$, so we obtain
\begin{equation}\label{1.6b}
\rank(C_{11}^{(2)})=\rank(B_{11}^{(2)})=r_e+r_b-n.
\end{equation}

Compute orthogonal matrices $U_2$ and $V_2$ such that
\[
U_2 B_{11}^{(2)}=\mat{c} 0 \\ B_{21} \rix, \quad C_{11}^{(2)} V_2 =\mat{cc} 0 & C_{12} \rix,
\]
where $B_{21}, C_{12}\in \R^{(r_e+r_b-n)\times (r_e+r_b-n)}$ are nonsingular.

Partition the matrices
\[ \bmat{ & n-r_e & r_e+r_b-n \cr
        n-r_e & E_{11} & E_{12} \cr
   r_e+r_b-n & E_{21} & E_{22} \cr}=U_2 E_{11}^{(1)}V_2,
\]
\[
U_2 B_{12}^{(2)}=\bmat{ &      \cr                                            n-r_e           & B_{12} \cr                                                     r_e+r_b-n     & B_{22} \cr}, \quad
C_{21}^{(2)}V_2=\bmat{ & n-r_e & r_e+r_b-n \cr
                & C_{21} & C_{22} \cr},
\]
and
\[
B_{32}=B_{22}^{(2)}, \quad C_{23}=C_{22}^{(2)}, \quad
    U=\mat{cc} U_2 & \\ & I \rix U_1, \quad
    V=V_1\mat{cc} V_2 & \\ & I \rix
\]
accordingly.

Using $E^TQ\geq 0$ and $C=B^TQ$  it follows that
\[
UQV=\bmat{  & n-r_b & r_e+r_b-n & n-r_e \cr
                     n-r_b   &  Q_{11} & Q_{12} & 0  \cr
                     r_e+r_b-n   & 0 & Q_{22} & 0  \cr
                     n-r_e   &  Q_{31} & Q_{32} & Q_{33}   \cr},
\]
as well as
\[
C_{12}=B_{21}^TQ_{22}, \ C_{23}=B_{32}^TQ_{33},
\]
which together with the nonsingularity of $C_{12}$, $C_{23}$, $B_{21}$  and $B_{32}$ implies that $Q_{22}$ and $Q_{33}$ are nonsingular.
In the following we show that $E_{11}$ is nonsingular by employing
that
\[
\mat{cc} Q_{11} & Q_{12}  \\ 0 & Q_{22} \rix^T \mat{cc} E_{11} & E_{12} \\ E_{21} & E_{22} \rix\geq 0.
\]
We may assume without loss of generality that
\[
Q_{22}=I,  \ Q_{11}=\mat{cc} I & 0 \\ 0 & 0 \rix, \
   Q_{12}=\mat{c} Z_1 \\ Z_2 \rix,
\]
and
\[
E_{11}=\mat{cc} E_{11}^1 & E_{11}^2 \\ E_{11}^3 & E_{11}^4\rix, \ E_{12}=\mat{c} E_{12}^1 \\ E_{12}^2 \rix, \
   E_{21}=\mat{cc} E_{21}^1 & E_{21}^2 \rix.
\]
Then
\[ \mat{cc} Q_{11} & Q_{12}  \\ 0 & Q_{22} \rix^T \mat{cc} E_{11} & E_{12} \\ E_{21} & E_{22} \rix
   =\mat{ccc} E_{11}^1 & E_{11}^2 & E_{12}^1 \\ 0 & 0 & 0 \\ \hat E_{21}^1 & \hat E_{21}^2 & \hat E_{22} \rix \geq 0,
\]
and thus,
\[
\rank(E_{11}^1)=\rank \mat{ccc} E_{11}^1 & E_{11}^2 & E_{12}^1 \rix, \ E_{11}^1 {\rm \ is \ nonsingular},
\]
and
\[ E_{11}^2=0, \
   \hat E_{21}^2=Z_1^TE_{11}^2+Z_2^T E_{11}^4 +E_{21}^2=Z_2^TE_{11}^4+E_{21}^2=0,
\]
i.e.,
\[ E_{11}^2=0, \ E_{21}^2=-Z_2^T E_{11}^4, \ \rank(E_{11}^4)=\rank\mat{c} E_{11}^4 \\ E_{12}^2 \rix, \ \mbox{and}\
    E_{11}^4 {\rm \ is \ nonsingular}. \]
Hence, $E_{11}$ is nonsingular.
\eproof

\subsection {Proof of Lemma~\ref{lemma-3}}

\proof Set
\begin{eqnarray*} X&=& \mat{ccc} I & 0 & -B_{12}B_{32}^{-1} \\ 0 & I & 0 \\ 0 & 0 & I \rix
\mat{ccc} I &     &       \\ \cU_{21} & \cU_{22} & \\ &  & I \rix U, \\
  Y&=&V \mat{ccc} I & \cV_{12} & \\
  & \cV_{22} &  \\   &  & I \rix
\mat{ccc} I & 0 & 0 \\ 0 & I & 0 \\ -C_{23}^{-1}C_{21} & 0 & I \rix .
\end{eqnarray*}
Then
\begin{eqnarray}
X E Y
 &=&\mat{ccc}            E_{11} & 0 & 0  \\
                              0 & \hat E_{22} & 0  \\
                                 0 & 0 & 0   \rix, \quad
  X A Y
 =\mat{ccc}            \tilde A_{11} & \tilde A_{12} & \tilde A_{13} \\
                              \tilde A_{21} & A_{22} & A_{23}  \\
                                \tilde A_{31} & A_{32} & A_{33}  \rix,
\nonumber\\
X  B W&=&\mat{cc}0         & 0 \\
\hat B_{21} & \hat B_{22} \\
0                & B_{32} \rix, \quad \quad
 W^TC Y
   =\mat{ccc}
          0 &  \hat C_{12} & 0 \\
          0 & \hat C_{22} & C_{23} \rix. \label{1.7}
\end{eqnarray}
Note that $E_{11}$ is nonsingular, so $\cU_{22}$ and $\cV_{22}$ are nonsingular, see \cite{ChuLM03}. Hence, $\hat B_{21}$ and $\hat C_{12}$ are nonsingular. It follows from
$E^TQ\geq 0$ and $C=B^TQ$ that
\begin{equation*} X^{-T} Q Y =\bmat{  & n-r_b & r_e+r_b-n & n-r_e \cr
                     n-r_b   &  \cQ_{11} & 0 & 0  \cr
                     r_e+r_b-n   & 0 & \cQ_{22} & 0  \cr
                     n-r_e   &         0 & \cQ_{32} & \cQ_{33}   \cr},
\end{equation*}
where $\cQ_{22}$ and $\cQ_{33}$ are nonsingular, and
\begin{equation*}
E_{11}^T\cQ_{11}\geq 0, \ \hat E_{22}^T\cQ_{22}\geq 0.
\end{equation*}
From $C=B^TQ$ and $A^TQ+Q^TA\leq 0$ it follows that
\begin{eqnarray}\label{1.10}
&& \mat{cc} \hat C_{12} & 0 \\ \hat C_{22} & C_{23} \rix=\mat{cc} \hat B_{21} & \hat B_{22} \\ 0 & B_{32} \rix^T \mat{cc} \cQ_{22} & 0 \\ \cQ_{32} & \cQ_{33} \rix, \\
&& \mat{cc} A_{22} & A_{23} \\ A_{32} & A_{33} \rix^T \mat{cc} \cQ_{22} & 0 \\ \cQ_{32} & \cQ_{33} \rix
  + \mat{cc} \cQ_{22} & 0 \\ \cQ_{32} & \cQ_{33} \rix^T  \mat{cc} A_{22} & A_{23} \\ A_{32} & A_{33} \rix \leq 0. \nonumber
\end{eqnarray}
This implies that
\[
Z \mat{cc} \cQ_{22} & 0 \\ \cQ_{32} & \cQ_{33} \rix \cZ=\bmat{ & \mu & r_b-\mu \cr
\mu                & \hat \cQ_{22} & 0 \cr
r_b-\mu           & \hat \cQ_{32} & \hat \cQ_{33} \cr}.
\]
%where $\hat \cQ_{22}$ and $\hat \cQ_{33}$ are nonsingular,
Then it follows that
\[
\cW^T \mat{cc} \hat C_{12} & 0 \\ \hat C_{22} & C_{23} \rix \cZ=(Z \mat{cc} \hat B_{21} & \hat B_{22} \\ 0 & B_{32} \rix \cW)^T
\mat{cc} \hat \cQ_{22} & 0 \\ \hat \cQ_{32} & \hat \cQ_{33} \rix
 =\mat{cc} \cC_{12} & 0 \\ \cC_{22} & \cC_{23} \rix,
\]
where $ \cC_{12}=\cB_{21}^T\hat \cQ_{22}$ and $\cC_{23}=\cB_{32}^T\hat \cQ_{33}$ are nonsingular and,
furthermore, $ \hat \cQ_{22}$ and $\hat \cQ_{33}$ are also nonsingular.
\eproof

\end{document}